\documentclass{mrlart3}
 \def\volno{xx}
  \def\yrno{2007}
  \def\issueno{x}
  \def\lpageno{100NN} 
\newtheorem{theorem}{Theorem}

\newtheorem{prop}[theorem]{Proposition}

\theoremstyle{remark}

\newcommand{\tit}{Spectrum of the Laplacian on manifolds with Spin(9) holonomy}

\renewcommand{\bar}{\overline}

\newcommand{\eps}{\varepsilon}
\newcommand{\ka}{K\"ahler }

\newcommand{\qk}{quaternionic \ka}

\newcommand{\nonum}{\nonumber }

\newcommand{\goto}{\rightarrow }

\newcommand{\re}{{\mathbb R}}

\newcommand{\bea}{\begin{eqnarray} }
\newcommand{\eea}{\end{eqnarray} }
\newcommand{\beay}{\begin{eqnarray*} }
\newcommand{\eeay}{\end{eqnarray*} }
\newcommand{\benu}{\begin{enumerate}}
\newcommand{\eenu}{\end{enumerate}}
\newcommand{\barr}{\begin{array}}
\newcommand{\earr}{\end{array}}
\newcommand{\grad }{\nabla}
\newcommand{\lap}{\triangle }
\newcommand{\pa }{\partial }
\newcommand{\sub }{\subseteq}

\newcommand{\lvl}{\mathcal {L}}

\newcommand{\la}{\langle}
\newcommand{\ra}{\rangle}

\newcommand{\cay}{\mathbb H_{\mathbb O}^2}
\newcommand{\oct}{\mathbb O}
\newcommand{\lt}{\left}
\newcommand{\rt}{\right}
\newcommand{\bay}{\begin{array}}
\newcommand{\eay}{\end{array}}
\newcommand{\ric}{\mbox{Ric}}
\newcommand{\beg}{\begin}
\newcommand{\hol}{\mbox{Hol($p$)}}

\title{\tit}
\author{Kwan-hang Lam}

\begin{document}
\maketitle \author

\begin{abstract} We consider noncompact complete manifolds with Spin(9)
holonomy and proved an one end result and a splitting type theorem
under different conditions on the bottom of the spectrum. We proved
that any harmonic functions with finite Dirichlet integral must be
Cayley-harmonic, which allowed us to conclude an one end result. In
the second part, we established a splitting type theorem by
utilizing the Busemann function.
\end{abstract}

\section*{Introduction}In \cite{klz}, the authors
proved the following \beg{thm}\cite{klz} Let $M$ be a complete Riemannian manifold with a
parallel $p$-form $\omega $. Assume that $f$ is a harmonic function satisfying
$$\int_{B_p(R)}|\grad f|^2=o(R^2)$$ as $R\goto \infty , $ then $f$ satisfies
\beay d*(df\wedge \omega)=0. \eeay\end{thm} Combining the above
theorem with the fact that a \qk manifold supports a global parallel
4-form $\omega $, the authors proved, by an explicit calculation
involving $\omega $, that a harmonic function with bounded Dirichlet
integral is quaternionic-harmonic. Utilizing the
quaternionic-harmonic condition they proved that, under an
assumption on the bottom of the spectrum $\lambda_1(M),$ such a
manifold must have exactly one infinite volume end. Since a manifold
with holonomy group Spin(9) supports a global parallel 8-form
$\Omega $, by a careful and detail study of $\Omega$, we proved that
any harmonic functions with bounded Dirichlet integral is
Cayley-harmonic. Similar to the work in \cite{klz}, with a suitable
lower bound assumption on $\lambda_1(M)$, an one infinite volume end
result has been established by utilizing the Cayley-harmonicity
condition. In the second part of this paper, we consider the case
that $\lambda_1(M)=121$ achieves its maximal value. By studying the
Busemann function $\beta $ on $M$ and using the results in
\cite{liwangpos} and \cite{liwangweight}, we proved that either $M$
has only one end or $M$ must splits as $\re \times N,$ where $N$ is
given by a level set of $\beta .$

\section{Cayley hyperbolic space}
We first give a brief introduction on the Cayley numbers $\oct$, and
a description of the sectional curvature of the Cayley hyperbolic
space $\cay.$ The material presented here is adopted from
\cite{brown-gray}, we refer the readers to there for further
details. The Cayley numbers $\oct$, is an 8-dimensional
non-associative division algebra over the real numbers which
satisfies the alternative law: $x(xy)=x^2y,\ (yx)x=yx^2$. It has a
multiplicative identity 1 and a positive definite bilinear form
$\la,\ra$ whose associated norm $||\cdot ||$ satisfies
$||ab||=||a||\cdot ||b||.$ Every element $a\in \oct$ can be written
as $a=\alpha 1+a_0$, where $\alpha $ is real and $\la a_0,1\ra=0.$
The conjugation map $a\mapsto a^*=\alpha 1-a_0$ is an
anti-automorphism, that is $(ab)^*=b^*a^*.$ Moreover, $aa^*=\la
a,a\ra 1$ and $\la a,b \ra=\la a^*,b^* \ra$. $\oct$ admits a
canonical basis $\{1,e_0,\cdots, e_6\}$ such that $\la e_i,e_j
\ra=\delta _{ij}$, $e_i^2=-1$, $e_ie_j+e_je_i=0$ for $i\not=j,$ and
$e_ie_{i+1}=e_{i+3},$ if $i$ is an integer mod 7. Obviously, we can
extend the positive bilinear form from $\oct$ to $\oct^2$ by
$$\la(a,b),(c,d)\ra=\la a,c\ra+\la b,d\ra,$$ where $a,b,c,d\in \oct.$ For any
point $x\in \cay,$ we make the following identification $T_x\cay
\simeq\oct^2.$ Let $M$ be a Riemannian manifold with metric tensor
$\la,\ra.$ Let $V$ be any tangent space to $M$. The curvature
operator of $M$ at $V$ is a map $$R:\Lambda ^2(V)\goto
\Lambda^2(V)\subseteq \mbox{Hom($V,V$)}$$ such that
$$R(x\wedge y)z+R(z\wedge x)y+R(y\wedge z)x=0. $$ The above two properties
implies $R$ is a symmetric linear operator, that is \beay \la
R(x\wedge y)z,w \ra = \la R(x\wedge y),z\wedge w \ra =\la R(z\wedge
w),x\wedge y \ra \eeay for any $x,y,z,w\in V.$ For any $x,y\in V$
linearly independent, the sectional curvature of the 2-plane spanned
by $x$ and $y$ is defined by
$$K_{x\wedge y}=\frac{\la R(x\wedge y),x\wedge y\ra}{||x\wedge y||^2}.$$
The sectional curvature $K_{(a,b)\wedge(c,d)}$ of the 2-plane
$(a,b)\wedge (c,d)$ of $\oct^2$ has the following properties:

\beg{enumerate}\item For any $a,b,c,d \in \oct $ with
$||(a,b)||=||(c,d)||=1$ and $\la (a,b),(c,d)\ra =0,$ we have \beay
K_{(a,b)\wedge(c,d)}&=&\alpha \big\{ ||a\wedge c||^2+||b\wedge
d||^2+\frac{1}{4}||a||^2||d||^2+\frac{1}{4}||b||^2||c||^2 \\
&+&\frac{1}{2}\la ab,cd\ra -\la ad,cb\ra\big\}\eeay

\item $$K_{(a,0)\wedge(b,0)}=\alpha \ \ \mbox{if}\ \ (a,0)\wedge(b,0)\not=0.$$

\item $$K_{(a,0)\wedge(0,b)}=\frac{\alpha}{4} \ \ \mbox{if}\ \ (a,0)\wedge(0,b)\not=0.$$

\item $$\frac{|\alpha|}{4}\leq |K_{(a,b)\wedge(c,d)}|\leq |\alpha |\ \ \mbox{if}\ \ (a,b)\wedge(c,d)\not=0.$$
\end{enumerate} In this article, we use the normalization that
$\alpha =-4$, hence the sectional curvature of $\cay$ is pinched
between $-4$ and $-1.$ Let $M$ be a complete noncompact Riemannian
manifold with holonomy group Spin(9). It was proved in
\cite{brown-gray} that a manifold with holonomy group Spin(9) must
be locally symmetric and its universal covering is either the Cayley
projective plane or the Cayley hyperbolic space $\cay.$ Since we are
considering noncompact manifolds, its universal covering is $\cay.$
We first compute the Laplacian of the distance function of $\cay $.
\begin{prop}\label{p:lap}Let $r(x)=r_p(x)$ be the distance function of $\cay$
from a fixed point $p$, then $$\lap r=14 \coth 2r+8\coth
r.$$\end{prop}
\begin{proof}
Let $\gamma:[0, L]\goto M$ be a normal geodesic from $p$ to $x$. Let
$e_1(t)=\gamma '(t)$ along $\gamma. $ Let $\{e_A\}_{A=2}^{16}$ be a
basis of $T_p\cay$ such that \bea \lt\{\bay{cc}
\label{e:cur1}{R_{1i1i}}= -4, &2\leq i\leq 8\\
R_{1\alpha 1\alpha }=-1, &9\leq \alpha \leq 16.\eay \rt.\eea We
extend $\{e_A\}$ to be a local frame along $\gamma (t)$, $\{\gamma
'(t)=e_1(t),e_2(t),\cdots, e_{16}(t)\}$ by parallel transporting
along $\gamma .$ Since $\cay $ is a symmetric space and thus locally
symmetric, we have
$$\frac{\partial }{\partial t}R_{1A1A}=R_{1A1A,1}=0, \ 2\leq A\leq 16,$$ hence
(\ref{e:cur1}) is valid along $\gamma .$ Let $X_A (t)=f_A(t)e_A(t)$
be the Jacobi field along $\gamma $ with $X_A(0)=0,\ X_A(L)=e_A(L).$
$f_A(t)$ satisfies the Jacobi equation \beay
&&\frac{d^2}{dt^2}f_A(t)-c_A^2f_A(t)=0\\
&&f_A(0)=0,\ f_A(p)=1, \ \ 2\leq A\leq 16.\eeay where $c_i=2 ,\
2\leq i\leq 8$ and $c_\alpha =1,\ 9\leq \alpha \leq 16.$ Solving the
above equation, we have \bea f_A(t)=\frac{\sinh (c_At)}{\sinh
(c_AL)},\ 2\leq A\leq 16.\eea Now, we can compute the Hessian of $r$
at $x$\beay
H(r)(e_A,e_A)&=&\int_0^L\left(\left|\frac{dX_A}{dt}\right|^2-\la
R(X_A,\gamma
')\gamma ',X_A\ra\right)dt\\
&=&\int_0^L\left(\left|\frac{df_A}{dt}\right|^2+c_A^2f^2\right)dt\\
&=&c_A\coth (c_AL).\eeay Therefore, we conclude that \beay \lap
r&=&\sum
_{A=2}^{16}H(r)(e_A,e_A)\\
&=&14\coth 2r+8\coth r,\eeay where we have used the fact that $
H(r)(e_1,e_1)=0.$
\end{proof}

\begin{theorem}\label{thm:lap} Let $M$ be a locally symmetric space with universal covering $\cay .$
Then $$\lambda _1(M)\leq 121,$$ and
$$\lap _Mr\leq 14\coth 2r+8\coth r,$$ in the sense of distribution. \end{theorem}
\begin{proof}Let $A(r), V(r)$ be the area and volume of the geodesic ball of radius $r$ of
$\cay$ respectively. By proposition \ref{p:lap}, we have $$
\frac{A'(r)}{A(r)}=14\coth 2r+8\coth r,$$ hence \bea
\label{e:volume}
V_M(p,r)&\leq &V(r)\\
\nonum &=&\int_0^rA(t)dt\\
\nonum &\leq &C \int_0^r(\sinh 2t)^7(\sinh t)^8dt\\
 \nonum&\leq &C_1e^{22r},\eea where $V_M(p,r)$ is volume of the geodesic ball with
radius $r$ centered at $p$ and for some constant $C_1$. On the other
hand, it was shown in \cite{liwangpos} that $$V_M(p,r)\geq C_2 \exp
\left({2\sqrt {\lambda _1(M)}r}\right),$$ for any manifolds with
positive spectrum. Combining the above inequality with
(\ref{e:volume}), we conclude that $\lambda_1(M)\leq 121.$  For the
second part, let $f(r)=14\coth 2r+8\coth r.$ By proposition
\ref{p:lap}, we have
$$\lap_Mr(x) =f(r(x)),$$  for any $x\in M\setminus
\mbox{Cut} (p),$ where Cut$(p)$ is the cut locus of $p$. For each
direction $\theta \in S_p(M),$ let $R(\theta ) =\sup
_{t>0}\{t:r_p(\exp _p(t\theta ))=t\}.$ Let $\phi \in C^\infty _0(M)$
be a non-negative smooth function with compact support, then \beay
\int_M\phi f(r)&=&\int _{S_p(M)}\int_0^{R(\theta
)}\phi f(r)\ J(\theta, r)\ drd\theta \\
&=&\int _{S_p(M)}\int_0^{R(\theta )}\phi  \ \frac{\pa J}{\pa r}\
drd\theta
\\
&=&-\int_M\frac{\pa \phi }{\pa r}+\int_{S_p(M)}\phi (\theta
,R(\theta
))J(\theta, R(\theta ))\ d\theta \\
&\geq &-\int_M\la \grad \phi ,\grad r\ra \\
&=&\int_Mr\lap \phi ,\eeay where the second equality follows from
the fact that $\lap r =\frac{\pa }{\pa r}(\log J),$ for all
$r<R(\theta)$ and the third equality follows from integration by
parts, $\phi \geq 0$ and $J(\theta ,0)=0$. Hence the second result
follows.
\end{proof}

Let us recall the definition of the Busemann function and some of
its properties. Let $M$ be a complete manifold and $\gamma
:[0,+\infty )\goto M$ be a geodesic ray. Let $\beta ^t_\gamma
(x)=t-r(\gamma (t),x),$ where $r(x,y)$ denotes the distance between
$x$ and $y$. Triangle inequality implies
$$|\beta ^t_\gamma (x)|=|r(\gamma (t),\gamma (0))-r(\gamma (t), x)|\leq r(\gamma
(0),x),$$ and $$\beta ^t_\gamma (x)-\beta ^s_\gamma (x)=t-s+r(\gamma
(s), x)-r(\gamma (t), x)\geq 0,$$ if $t>s.$ Hence $\{\beta ^t_\gamma
\}_{t\geq 0}$ is uniformly bounded on compact subsets of $M$ and
nondecreasing, it converges uniformly on any compact subsets of $M$.
The Busemann function with respect to a geodesic ray $\gamma $ is
defined as $$\beta (x)=\lim _{t\goto +\infty }\beta _\gamma^t(x).$$
The following lemma is well-known and the proof here is adopted from
\cite{liwangsymm}.

\beg{lem} \label{lem:busemann}$$\lt|\grad \beta \rt|=1,$$ almost
everywhere.
\end{lem} \beg{proof}Triangle inequality implies
$$|\beta ^t_\gamma (x)-\beta ^t_\gamma (y)|\leq r(x,y),$$ which implies $\beta
$ is Lipschitz with Lipschitz constant 1. For any point $x\in M$, we
consider a normal geodesic $\tau _t$ joining from $x=\tau _t(0)$ to
$\gamma (t).$ Since the unit sphere is compact,
$\{\tau_t'(0)\}_{t>0}$ has a limit point $v\in T_xM.$ The sequence
$\tau _t$ converges to a geodesic ray $\tau $ with $\tau (0)=x$ and
$\tau '(0)=v.$ Hence, if we let $s,\eps >0$, if $t$ is sufficiently
large, we have $r(\tau_t(s),\tau(s))<\eps .$ Again, triangle
inequality implies \beay \beta (\tau (s))-\beta (\tau (0))&=&\lim
_{t\goto \infty }(r(\tau (0),\gamma
(t))-r(\tau (s),\gamma (t)))\\
&=&\lim_{t\goto \infty }(r(\tau (0),\gamma (t))-r(\tau_t(s),\gamma
(t))+
r(\tau_t(s),\gamma (t))-r(\tau (s),\gamma (t)))\\
&\geq &\lim_{t\goto \infty }(r(\tau (0),\gamma
(t))-r(\tau_t(s),\gamma
(t))-r(\tau_t (s),\tau(s)))\\
&\geq &\lim_{t\goto \infty }(r(\tau (0),\gamma
(t))-r(\tau_t(s),\gamma
(t)))-\eps \\
&\geq& s-\eps ,\eeay thus \bea \label{e:bus}|\beta (\tau (s))-\beta
(\tau (0))|\geq s.\eea The result follows by combining the above
inequality with the fact that $\beta $ is a Lipschitz function with
Lipschitz constant 1.\end{proof}

\section{Manifolds with a parallel form}

Let us first recall the Hodge star operator * and some of its basic
properties. Let $V^n$ be a $n$-dimensional oriented real inner
product space, we have the Hodge star operator $$*:\wedge ^pV\goto
\wedge ^{n-p}V,$$ for any $\theta \in \wedge ^1V, v\in V,$ exterior
multiplication and interior product operators \beay
\eps (\theta ):&&\wedge ^pV\goto \wedge ^{p+1}V\\
l(v):&&\wedge ^pV\goto \wedge ^{p-1}V,\eeay where $\eps
(\theta)\omega =\theta \wedge \omega$ and $\lt(l(v)\omega \rt)(\cdot
)=\omega (v,\cdot)$ for any $\omega \in \wedge^pV.$ Let
$\theta,\theta '\in \wedge ^1V$ and $v,v'\in V$ be the dual of
$\theta $ and $\theta '$ respectively with respect the inner product
of $V$. For any $\eta \in \wedge ^pV,$ we have the following basic
properties \beg{enumerate}
\item $**\eta =(-1)^{p(n-p)}\eta$
\item $*\eps (\theta )\eta =(-1)^pl(v)*\eta $
\item$\eps (\theta)*\eta =(-1)^{p-1}*l(v)\eta$
\item$ *\eps (\theta)*\eta =(-1)^{(p-1)(n-p)}l(v)\eta $
\item$l(v)\eps (\theta ')\eta +\eps (\theta )l(v')\eta =0,$ where $v\perp v'$
\item$l(v)\eps (\theta)\eta +\eps (\theta )l(v)\eta =\eta $
\end{enumerate}

The following theorem is an over-determined system of equations
satisfied by harmonic functions and generalized Corlette's argument
to harmonic functions with finite Dirichlet integral on a complete
manifold with a parallel $p$-form. This kind of result was first
proved by Siu \cite{siu} for harmonic maps in his proof of the
rigidity theorem for \ka manifolds. Corlette \cite{cor} gave a more
systematic approach for harmonic maps with finite energy from a
finite volume quaternionic hyperbolic space or Cayley hyperbolic
plane to a manifold with nonpositive curvature. In \cite{likahler},
the author generalized Siu's argument to harmonic functions with
finite Dirichlet integral on \ka manifolds.

\beg{theorem}(\cite{klz})\label{thm:klz} Let $M$ be a complete
Riemannian manifold with a parallel $p$-form $\omega $. Assume that
$f$ is a harmonic function satisfying
$$\int_{B_p(R)}|\grad f|^2=o(R^2)$$ as $R\goto \infty , $ then $f$ satisfies
\beay d*(df\wedge \omega)=0. \eeay\end{theorem}

By taking a careful and closer look at the nature of the proof of
the above theorem, we found out that the proof not only works for
harmonic functions with finite Dirichlet integral but also $L^2$
harmonic 1-form. The key ingredient is that any $L^2$ harmonic
1-form is both closed and co-closed. We have the following:

\beg{theorem}\label{thm:cor1} Let $M$ be a complete Riemannian
manifold with a parallel $p$-form $\omega $. Assume that $\alpha $
is a $L^2$ harmonic 1-form, that is $\lap \alpha =0$ and
$$\int_{M}|\alpha |^2<+\infty .$$  Then $\alpha $ satisfies
\beay d*(\alpha \wedge \omega)=0. \eeay\end{theorem} \beg{proof} We
first show that \bea\label{e:para1} *d*(\alpha \wedge \omega
)=(-1)^{n-1}d*(\alpha \wedge
*\omega).\eea For any $x\in M$, we choose a local orthonormal frame
$\{e_i\}_{i=1}^n$ such that $\grad _{e_i}e_j(x)=0.$ Let $\{\theta
^i\}_{i=1}^n$ be the coframe. For any $p$-form $\omega $, we have
$$d\omega =\eps (\theta^i)\grad_{e_i}\omega$$ at $x$ and $\omega $
is parallel if and only if $\grad _{e_i}\omega =0, \forall i.$ Let
$\alpha =\sum _{i=1}^na_i\theta ^i,$ and hence $\bar \alpha
=\sum_{i=1}^na_ie_i$ is the dual of $\alpha .$ We use the notation
$d\alpha =\sum_{i,j=1}^na_{i,j}\theta ^j\wedge \theta ^i,$ where
$a_{i,j}=\grad _{e_j}a_i.$ Since $\alpha $ is $L^2$ harmonic, it is
both closed and co-closed, which are equivalent to the conditions
that $a_{i,j}=a_{j,i}$ and $\sum_{i=1}^na_{i,i}=0.$ The following
calculations are all evaluated at $x$. \bea \label{e:para2}d*(\alpha
\wedge
*\omega )&=&d*\eps (\alpha )*\omega \\
\nonum&=&(-1)^{(p-1)(n-p)}d[l(\bar \alpha )\omega ]\\
\nonum&=&(-1)^{(p-1)(n-p)}\sum_{i=1}^n\eps (\theta
^i)\grad_{e_i}(l(\bar \alpha
)\omega )\\
\nonum&=&(-1)^{(p-1)(n-p)}\sum_{i,j=1}^n\eps (\theta
^i)a_{j,i}(l(e_j)\omega ),\eea where the third equality follows from
$\grad_{e_i}e_j(x)=0$ and the last equality follows from $\grad
\omega =0$. On the other hand, \bea \label{e:para3}*d*(\alpha \wedge
\omega)&=&*d*\eps (\alpha
)\omega \\
\nonum&=& *\sum_{i=1}^n\eps (\theta ^i)\grad_{e_i}\lt(*\eps
\lt(\sum_{j=1}^na_j\theta ^j\rt)\omega \rt) \\
\nonum&=&*\lt(\sum_{i,j=1}^na_{j,i}\eps (\theta ^i)*\lt(\eps (\theta
^j)\omega
\rt)\rt)\\
\nonum&=&(-1)^{p(n-p-1)}\sum_{i,j=1}^na_{i,j}l(e_i)\eps (\theta
^j)\omega
\\
\nonum &=&(-1)^{p(n-p-1)}\lt(\sum_{i=1}^na_{i,i}l(e_i)\eps (\theta
^i)\omega
+\sum_{i\not=j}^na_{i,j}l(e_i)\eps (\theta ^j)\omega\rt)\\
\nonum&=&(-1)^{p(n-p-1)}\Big(\sum_{i=1}^na_{i,i}[\omega -\eps
(\theta
^i)l(e_i)\omega] \\
\nonum&&-\sum_{i\not=j}^na_{i,j}\eps (\theta ^j)l(e_i)\omega\Big)\\
\nonum&=&(-1)^{p(n-p-1)+1}\sum_{i,j=1}^na_{i,j}\eps (\theta
^i)(l(e_j)\omega ),\eea where the last equality follows from
$a_{i,j}=a_{j,i}$ and $\sum _{i=1}^na_{i,i}=0.$ (\ref{e:para1}) now
follows from (\ref{e:para2}) and
(\ref{e:para3}). Let \beay \phi (x)=\lt\{\bay{clc}1&\mbox{on }&B_p(R)\\
0&\mbox{on }&M\setminus B_p(2R)\eay\rt.\eeay such that $|\grad \phi
|\leq C_1R^{-1}.$ Consider \bea\label{e:para4} \ \ \ \ \ \
\int_M\phi ^2|d*(\alpha \wedge \omega )|^2&=&\lt|\int_M\phi
^2d*(\alpha \wedge \omega )\wedge
*d*(\alpha \wedge \omega
)\rt|\\
\nonum&=&\lt|\int_M\phi ^2d*(\alpha \wedge \omega )\wedge d*(\alpha
\wedge
*\omega
)\rt|\\
\nonum&=&\lt|\int_Md\phi ^2\wedge *(\alpha \wedge \omega )\wedge
d*(\alpha \wedge
*\omega )\rt|\\
\nonum&\leq &2\lt(\int_M|d\phi|^2|*(\alpha \wedge
\omega)|\rt)^{1/2}\lt(\int_M\phi
^2|d*(\alpha \wedge *\omega )|^2\rt)^{1/2}\\
\nonum&=&2\lt(\int_M|d\phi|^2|*(\alpha \wedge
\omega)|\rt)^{1/2}\lt(\int_M\phi ^2|d*(\alpha \wedge \omega
)|^2\rt)^{1/2},\eea where the second and the last equality follows
from (\ref{e:para1}), the third equality follows from integration by
parts and the fact that $d^2=0.$ $\omega $ is parallel implies
$$|*(\alpha \wedge \omega )|\leq C_2|\alpha |,$$ for some constant $C_2.$
Combining the above with (\ref{e:para4}), we have

$$\int_{B_p(R)}|d*(\alpha
\wedge \omega )|^2\leq C_3R^{-2}\int_{B_p(2R)}|\alpha |^2.$$ Let
$R\goto +\infty,$ the result follows from the assumption that
$\alpha$ is $L^2$ integrable. \end{proof}

\section{Some vanishing theorems}
The following lemma is useful in proving vanishing theorems
\beg{cor}\cite{lamwpi}\label{cor:vanish} Let $b
>-1$. Assume that $h$ is $L^2$ integrable and satisfies differential
inequality
$$\lap h\geq -ah+b\frac{|\grad h|^2}{h},$$ for some constant $a$. If $\lambda_1(M)>0$
and the Ricci curvature satisfies
$$\ric_M\geq -(b+1)\lambda_1(M)+\delta,$$ for some $\delta
>0$. Then $h\equiv 0.$ \end{cor}

Combining theorem \ref{thm:cor1} with corollary \ref{cor:vanish}, a sharper form of vanishing
theorems (\cite{liwangpos}, \cite{lamwpi}) for manifolds with a parallel $p$-form can now be
established: \beg{theorem}\label{thm:kahler vanish} Let $M^{2n}$ be a $2n$ real dimensional
\ka manifold with $\lambda_1(M)>0.$ Assume the Ricci curvature of $M$ satisfies
$$\ric_M\geq -2\lambda_1(M)+\delta ,$$ for some $\delta >0.$ Then
$H^1(L^2(M))=0.$\end{theorem} \beg{proof}Let $\omega \in
H^1(L^2(M))$ and $h=|\omega |.$ We claim that $h$ satisfies the
Bochner formula of the following form
$$\lap h\geq \frac{\ric_M(\omega, \omega )}{h}+\frac{|\grad h|^2}{h}.$$
Applying corollary \ref{cor:vanish} with $b=1,$ the result follows.
To prove the claim, we let $\{e_i\}_{i=1}^{2n}=\{\bar e_1,\cdots
,\bar e_n,I\bar e_1,\cdots I\bar e_n\} $ be a local orthonormal
frame, where $I$ is the complex structure and $\{\theta
^i\}_{i=1}^{2n}=\{\bar \theta ^1,\cdots, \bar \theta ^n,I\bar \theta
^1,\cdots,I\bar \theta^n\}$ be the orthonormal coframe. The \ka
form, satisfying $\Omega (X,Y)=g(X,IY),$ is then given by \beay
\Omega &=&-\sum_{i=1}^n\bar \theta ^i\wedge
I\bar \theta ^i\\
&=&-\sum_{i=1}^n\theta ^i\wedge \theta ^{n+i}.\eeay With the above
notations, we can write $\omega =\sum_{i=1}^{2n}a_i\theta ^i.$ By
theorem \ref{thm:cor1},
$$d*(\omega \wedge \Omega)=0,$$ which is equivalent to, by (\ref{e:para3})
\bea \label{e:kahler vanish}\sum_{i,j=1}^{2n}a_{ij}\theta ^i\wedge
l(e_j)\Omega =0,\eea where we have used the notation
$a_{ij}=a_{i,j}.$ Since \beay l(e_j)\Omega=\left\{\bay{cc}-\theta
^{j+n}&1\leq j\leq
n\\
\theta ^{j-n}&n+1\leq j\leq 2n\eay\right.,\eeay hence (\ref{e:kahler
vanish}) becomes \beay \sum_{i=1}^{2n}\lt( -\sum_{j=1}^na_{ij}\theta
^i\wedge \theta^ {j+n}+\sum_{j=1}^na_{i,{j+n}}\theta ^i\wedge
\theta^ {j}\rt)=0.\eeay The coefficient of $\theta ^i\wedge \theta
^{i+n}$ of the above equation is zero and thus we conclude that \bea
\label{e:kahler 1}a_{ii}+a_{i+n ,i+n}=0,\eea for any $1\leq i\leq
n$. Now we go back to study the form $\omega .$ Let
$\{e_i\}_{i=1}^{2n}$ as described above with $e_1$ such that $\omega
(e_1)=|\omega |$ and $\omega (e_j)=0$ for any $j\not=1$ at a fixed
point $p.$
\beay |\grad \theta |^2&=&\sum_{i,j=1}^{2n}a_{ij}^2\\
&\geq &a_{11}^2+a_{n+1,n+1}^2+2\sum_{j=2}^{2n}a_{1j}^2\\
&=&2\lt(a_{11}^2+\sum_{j=2}^{2n}a_{1j}^2\rt)\\
&=&2|\grad h|^2\eeay at $p,$ where the third equality follows from
(\ref{e:kahler 1}). Combining the above inequality with the Bochner
formula gives us \beay \frac{1}{2}\lap
(h^2)&=&\ric(\omega, \omega )+|\grad \theta |^2\\
&\geq &2|\grad h|^2+\ric (\omega, \omega ).\eeay Hence
$$\lap h\geq \frac{\ric_M(\omega, \omega )}{h}+\frac{|\grad h|^2}{h},$$ and
the claim is justified.

\end{proof}

\beg{theorem}\label{thm:qk vanish} Let $M^{4n}$ be a $4n$ dimensional \qk manifold. Assume
that $\lambda_1(M)>0$ and the Ricci curvature of $M$ satisfies
$$\ric_M\geq -\frac{4}{3}\lambda_1(M)+\delta ,$$ for some $\delta >0.$ Then
$H^1(L^2(M))=0.$\end{theorem} \beg{proof} We follow the notations in
\cite{klz}. $M$ has a rank 3 vector bundle $V\sub End(TM)$
satisfying \beg{enumerate}\item In a local coordinate neighborhood,
there exists a local basis $\{I,J,K\}$ of $V$ such that \beay
I^2=J^2=K^2=-1\\
IJ=-JI=K\\
JK=-KJ=I\\
KI=-IK=J\eeay and \beay g(X,Y)=g(IX,IY)=g(JX,JY)=g(KX,KY),\eeay for
any $X,Y\in TM.$ \item If $\phi \in \Gamma (V),$ then $\grad _X\phi
\in \Gamma (V)$ for any $X\in TM.$
\end{enumerate} We define following two forms \beay \omega_1(X,Y)&=&g(X,IY)\\
\omega_2(X,Y)&=&g(X,JY)\\
\omega_3(X,Y)&=&g(X,KY).\eeay The parallel 4-form of $M$ is then
given by \beay \Omega =\omega_1\wedge \omega_1+\omega_2\wedge
\omega_2+\omega_3\wedge \omega_3.\eeay Let
$$\{e_i\}_{i=1}^{4n}=\{\bar e_1,\cdots, \bar e_n,I\bar e_1,\cdots,
I\bar e_n,J\bar e_1,\cdots, J\bar e_n,K\bar e_1,\cdots, K\bar e_n\}$$ be a local orthonormal
frame and $$\{\omega ^i\}_{i=1}^{4n}=\{\bar \theta^1,\cdots, \bar \theta^n,I\bar
\theta^1,\cdots, I\bar \theta^n,J\bar \theta^1,\cdots, J\bar\theta^n,K\bar \theta^1,\cdots,
K\bar \theta^n\}$$ be the orthonormal coframe. Let $\omega =\sum_{i=1}^{4n}a_i\omega^i\in
H^1(L^2(M)).$ Using the above formula of $\Omega $ and calculate as in theorem \ref{thm:kahler
vanish} (or see \cite{klz}), we have \bea\label{e:qk 1}
a_{ii}+a_{i+n,i+n}+a_{i+2n,i+2n}+a_{i+3n,i+3n}=0,\eea for any $1\leq i\leq n.$ We now proceed
as in theorem \ref{thm:kahler vanish}. Let $h=|\omega |.$ It is not difficult to see that $h$
satisfies the Bochner formula of the following form
$$\lap h\geq \frac{\ric_M(\omega, \omega )}{h}+\frac{1}{3}\frac{|\grad h|^2}{h}.$$
Applying corollary \ref{cor:vanish} with $b=\frac{1}{3},$ the result follows. Indeed, let
$\{e_i\}_{i=1}^{4n}$ as described above with $e_1$ such that $\omega (e_1)=|\omega |$ and
$\omega (e_j)=0$ for any $j\not=1$ at a point $p.$ We compute
\beay |\grad \theta |^2&=&\sum_{i,j=1}^{4n}a_{ij}^2\\
&\geq &a_{11}^2+a_{1+n,1+n}^2+a_{1+2n,1+2n}^2+a_{1+3n,1+3n}^2+2\sum_{j=2}^{2n}a_{1j}^2\\
&\geq &a_{11}^2+\frac{1}{3}(a_{1+n,1+n}+a_{1+2n,1+2n}+a_{1+3n,1+3n})^2+2\sum_{j=2}^{2n}a_{1j}^2\\
&=&\frac{4}{3}a_{11}^2+2\lt(\sum_{j=2}^{2n}a_{1j}^2\rt)\\
&=&\frac{4}{3}\lt(a_{11}^2+\sum_{j=2}^{2n}a_{1j}^2\rt)\\
&=&\frac{4}{3}|\grad h|^2\eeay at $p,$ where the third inequality
and the fourth equality follow from Schwarz's inequality and from
(\ref{e:qk 1}) respectively. Combining the above inequality with the
Bochner formula gives us \beay \frac{1}{2}\lap
(h^2)&=&\ric(\omega, \omega )+|\grad \theta |^2\\
&\geq &\frac{4}{3}|\grad h|^2+\ric (\omega, \omega ).\eeay Hence
$$\lap h\geq \frac{\ric_M(\omega, \omega )}{h}+\frac{1}{3}\frac{|\grad h|^2}{h}.$$
\end{proof}

\section{Holonomy and Spin(9) invariant}

We give a very brief introduction and list some basic principles about the holonomy group of a
Riemannian manifold. We refer the readers to \cite{besse} and the references therein for
further details. Most of the following introductory material are adopted from there. Let $p\in
M$ and $\gamma:[0,l]\goto M $ be a $C^1$-piecewise closed curve with $\gamma (0)=\gamma
(l)=p$. Let $\tau (\gamma ):T_pM\goto T_pM$ be the parallel transport along $\gamma .$ Since
parallel transport preserves inner product, $\tau (\gamma )$ is an element of $O(T_pM)$, the
orthogonal group of $T_pM$. Since the inverse of a curve $\gamma ^{-1}$ and the composition of
two curves $\gamma \cup \sigma $ satisfy $\tau (\gamma ^{-1})=(\tau (\gamma ))^{-1}$ and $\tau
(\gamma \cup \sigma )=\tau (\gamma )\circ \tau (\sigma )$. We can have the following
definition:

\beg{defi} The holonomy group (or the holonomy representation of $M$
at $p$) of a Riemannian manifold $(M,g)$ at $p$ is defined by
$$\hol=\{\tau (\gamma ): \gamma\in C^1\mbox{-piecewise closed curves of $M$ based at $p$}\},$$ the subgroup of the orthogonal
group $O(T_pM)$. \end{defi} On $M$, let us consider a tensor field
$\alpha $. If $\alpha $ is invariant by parallel transport, that is,
for any $p,q\in M$ and any curve $\gamma$ from $p$ to $q$, we have
$$\tau ^*(\gamma )(\alpha (p))=\alpha (q),$$ where $\tau ^*(\gamma
)$ is the tensorial extension of the parallel transport $\tau
(\gamma )$ along $\gamma.$ By the above definition, $\alpha (p)$ at
$T_pM$ is hence invariant by the tensorial extension of the holonomy
representation $\hol\subseteq O(T_pM).$ Conversely, given any tensor
on $T_pM$, if $\alpha _0$ is invariant under the tensorial extension
of $\hol$, we can construct a tensor field $\alpha $ on $M$ by the
formula $\tau ^*(\gamma )(\alpha (p))=\alpha (q).$ Since $\alpha _0$
is invariant under the tensorial extension of $\hol$, the above
definition is independent of the choice of the curve $\gamma $ and
thus it is well-defined. Clearly, $\alpha (p)=\alpha_0$. By the
above discussion, we have established a fundamental principle of
holonomy group. \beg{prop}\label{p:holonomy principle} Let $M$ be a
Riemannian manifold and we consider a fixed type $(r,s)$ tensors on
$M$. Then the following three properties are
equivalent:\beg{enumerate}

\item There
exists a tensor field of type $(r,s)$ which is invariant by parallel
transport

\item There exists $p\in M$ and a tensor $\alpha_0$ of type $(r,s)$ which is
invariant by the tensorial extension of type $(r,s)$ of the holonomy
representation $\hol.$

\item There exists a tensor field $\alpha $ of type $(r,s)$ which has zero
covariant derivative.
\end{enumerate}\end{prop}
\beg{proof} We have already established the equivalency of the first
two statements in the discussion above. For the last statement, it
can be seen easily via the formula $$(D\alpha )(X_1,\cdots
,X_s;X)=D_X(\alpha (X_1,\cdots, X_s))-\sum_{i=1}^s\alpha
(X_1,\cdots, D_XX_i,\cdots, X_s).$$ For any curve $\gamma$, let
$X_1,\cdots,X_s$ be vector fields parallel along $\gamma $ and
$X=\gamma '.$ Hence, the above equation becomes $$(D\alpha
)(X_1,\cdots ,X_s;X)=D_X(\alpha (X_1,\cdots, X_s)).$$ Therefore,
$D\alpha =0$ is equivalent to $D_X(\alpha (X_1,\cdots, X_s))$, which
implies $\alpha (X_1,\cdots, X_s)$ is constant along $\gamma .$
Conversely, for any tangent vector $X(p),$ we can choose a curve
$\gamma $ such that $\gamma '=X(p).$
\end{proof}

Let $M$ be a manifold with holonomy group Spin(9). We are now ready
to describe the parallel 8-form of $M.$ The parallel 8-form of
$\cay$ has been obtained by Brown and Gray in \cite{brown-gray}.
However, it is not easy to read off its properties for further
applications because their 8-forms are defined via integration. In
\cite{abe}, the authors defined an 8-form $\Omega $  and showed that
it is Spin(9) invariant. In \cite{klz} the authors used the explicit
formula of the parallel 4-form of a quaternionic \ka manifold and
proved that any harmonic function with finite Dirichlet integral is
quaternionic-harmonic. Similarly, we will combine the explicit
formula of $\Omega$ in \cite{abe} with a result in \cite{klz} to
conclude that any harmonic function with finite Dirichlet integral
is Cayley-harmonic. We now give a brief description of the Spin(9)
invariant 8-form $\Omega$ and we will follow the notations in
\cite{abe}. For any point $p\in \cay$, we identify the tangent space
at $p$ to the ordered pair of Cayley numbers,
$T_p(\cay)=\oct^2=\{(x,y):x,y\in \oct\}$. Let $\bar e_0=1,\bar
e_1,\cdots \bar e_7$ be a basis of $\oct$ as in \cite{yok}. For any
$x\in \oct,$ we let $x^{(2)}=(x,0)$ and $x^{(3)}=(0,x).$ Let
$\{v_i\}_{i=0}^7$ be the dual 1-forms of $\{\bar
e_i^{(2)}\}_{i=0}^7$ and $\{w_i\}_{i=0}^7$ be the dual 1-forms of
$\{\bar e_i^{(3)}\}_{i=0}^7$. Equivalently, we have \beay
v_i(\bar e_j^{(2)})&=&\delta _{ij}, \ \ v_i(\bar e_j^{(3)})=0\\
w_i(\bar e_j^{(2)})&=&0, \ \ w_i(\bar e_j^{(3)})=\delta _{ij},\eeay
for any $0\leq i,j\leq 7.$ Let $e_i=\bar e_{i-1}^{(2)}$ for $  1\leq
i\leq 8$ and $e_j=\bar e_{j-9}^{(3)}$ for $9\leq j\leq 16$ so that
$\{e_i\}_{i=1}^{16}$ is an orthonormal basis of $T_p\cay.$
$$\omega_{ij}=v_{\sigma (i)}\wedge v_{\sigma(j)},\ \ \ \eta_{ij}=w_{\tau
(i)}\wedge w_{\tau(j)},$$ for some functions $\sigma, \tau $ which
are given in \cite{abe}. For our purpose, we do not need to know the
explicit forms of $\sigma ,\tau$ and so we ignore it here for the
sake of simplicity. Now we are ready to write down the formula of
$\Omega.$ \beg{theorem}\cite{abe} With the above notations,
\label{thm:parallel form}
$$\Omega=(-v_0\wedge \cdots \wedge v_7+w_0\wedge \cdots \wedge w_7)+F(\omega_{ij},\eta_{kl})$$
is Spin(9) invariant, where $F$ is a linear combinations of 8-forms,
each of which is wedge products of some combinations of
$\omega_{ij},\eta_{kl}$.
\end{theorem} We would like to point out that $F$ was given explicitly in
\cite{abe}. However, the above simplified form of $\Omega $ is enough for our application.
\beg{theorem}\label{thm:cayley harmonic}Let $M$ be a manifold with holonomy group Spin(9).
Assume that $f$ is a harmonic function satisfying
$$\int_{B_p(R)}|\grad f|^2=o(R^2),$$ as $R\goto \infty.$ Then with the above notations, we have $$\sum _{i=1}^8f_{ii}=0,$$
where $f_{ij}=Hess(f)(e_i,e_j).$\end{theorem} \beg{proof} Fix $x\in
M$ and let $\{e_i\}_{i=1}^{16}$ be the orthonormal frame of $T_xM$
in the above discussion. By the above construction, let $\{\theta^i
\}_{i=1}^{16}=\{v_0,\cdots, v_7, w_0,\cdots, w_7\}$ be the
orthonormal coframe. By theorem \ref{thm:parallel form},
$$\Omega=(-v_0\wedge \cdots \wedge v_7+w_0\wedge \cdots \wedge
w_7)+F(\omega_{ij},\eta_{kl})$$ is Spin(9) invariant. Since $M$ has
holonomy group Spin(9), by proposition \ref{p:holonomy principle},
$\Omega$ can be extended to be a parallel form on $M$, which we
still denote it by $\Omega.$ By theorem \ref{thm:klz}, we have
$$d*(df\wedge \Omega)=0.$$ From (\ref{e:para3}), by replacing $a_{i,j}$ by $f_{ij}$, the above equation is
equivalent to $$\sum_{i,j=1}^{16}f_{ij}\eps (\theta
^i)\lt(l(e_j)\Omega\rt)=0.$$ Evaluate the above equation at $x$, we
claim that the only terms contain $v_0\wedge \cdots\wedge v_7$ are
the following \beay \sum_{i=1}^{8}f_{ii}\eps (\theta
^i)\lt(l(e_i)(-v_0\wedge \cdots \wedge
v_7)\rt)&=&-\sum_{i=1}^8f_{ii}v_0\wedge \cdots \wedge v_7.\eeay
Since the coefficient of $v_0\wedge \cdots \wedge v_7$ of
$d*(df\wedge \Omega)$ is zero, we conclude that
$$\sum_{i=1}^8f_{ii}=0,$$ at $x$. To prove the claim, since $l(e_j)F(\omega_{ab},\eta_{cd})$ kills off a
$v_{j-1}$ term if $1\leq j\leq 8 $ or a $w_{j-9}$ term $9\leq j\leq
16$ of $F(\omega_{ab},\eta_{cd})$. On the other hand, when $\eps
(\theta ^i)$ acts on $l(e_j)F(\omega_{ab},\eta_{cd})$, it adds a
$v_{i-1}$ term if $1\leq i\leq 8 $ or a $w_{i-9}$ term $9\leq i\leq
16$ to $l(e_j)F(\omega_{ab},\eta_{cd})$. Since
$$\omega_{ab}=v_{\sigma (a)}\wedge v_{\sigma (b)}$$ and
$$\eta_{ab}=w_{\tau (a)}\wedge w_{\tau (b)},$$ by the above discussion, for any $1\leq i,j\leq 16,$
  $\eps(\theta^i)\lt(l(e_j)F(\omega_{ab},\eta_{cd})\rt)$ does not contain any
terms of the form $v_0\wedge\cdots\wedge v_7$ and
$w_0\wedge\cdots\wedge w_7$. This proved the claim and the result
follows.
\end{proof}


\section{Manifolds with positive spectrum}
We will summarize some useful properties of manifolds with positive spectrum. We refer the
readers to \cite{liwangpos} for a more detailed description on this subject. Let $M$ be a
manifold with positive spectrum $\lambda_1 (M)>0.$ By the variational principle, it is
equivalent to the following condition:
$$\lambda _1(M)\int _M\phi ^2\leq \int_M|\grad \phi |^2,$$ for any compactly
supported smooth function $\phi \in C^\infty_c(M).$ Since
$\lambda_1(M)>0,$ $M$ must be nonparabolic and it implies $M$ must
have at least one nonparabolic end. $\lambda_1 (M)>0$ also implies
an end $E$ of $M$ is nonparabolic if and only if it has infinite
volume. Assume that $M$ has at least two infinite volume ends,
$E_1,E_2.$ Let $B_p(R)$ be the geodesic ball with radius $R$
centered at $p$. We write $B(R)=B_p(R)$ when there is no ambiguity.
We construct a sequence of harmonic functions $\{f_R\}$ by solving
the
following equation \beay \bay{rcc}\lap f_R=0 &\mbox{on}& B(R)\\
f_R=1&\mbox{on}&\pa B(R)\cap E_1\\
f_R=0&\mbox{on}&\pa B(R)\setminus E_1\eay.\eeay By the theory of \cite{litamstructure},
$\{f_R\}$ converges (by passing to a subsequence if necessary) to a nonconstant harmonic
function $f$ with finite Dirichlet integral on $M$ as $R\goto +\infty.$ Maximum principle
implies that $0\leq f\leq 1$. By the construction, it is clear that $\sup _Mf=\sup _{E_1}f=1$
and $\inf_Mf=\inf_{E_2}f=0$. We will need the following lemmas:

\beg{lem}\cite{liwangpos}\label{lem:est1-pos} With the above notations, $f$ as constructed
above. Then \beg{enumerate}\item \beay \int_{E_1(R+1)\setminus E_1(R)}(1-f)^2&\leq &C\exp
(-2\sqrt {\lambda_1
(M)}R)\\
\int_{E(R+1)\setminus E(R)}f^2&\leq &C\exp (-2\sqrt {\lambda_1
(M)}R)\eeay for some constant $C$ depends on $f,\ \lambda_1(M)$ and
the dimension of $M,$ where $E$ is any other end different from
$E_1.$
\item $$\int_{E(R+1)\setminus E(R)}|\grad f|^2\leq C\exp (-2\sqrt{\lambda_1 (M)}R),$$ for $R$ sufficiently
large, where $E$ is any end of $M.$
\end{enumerate}
\end{lem}

\beg{lem}\cite{liwangweight}\label{lem:est-weight} For the function $f$ constructed above, let
$\inf f<a<b<\sup f$,
$$l(t)=\{x\in M:f(x)=t\}$$ and $$\lvl(a,b)=\{x\in M:a<f(x)<b\}.$$ Then $$\int_{\lvl(a,b)}|\grad f|^2=(b-a)\int_{l(b)}|\grad f|$$ and
$$\int_{l(b)}|\grad f|=\int_{l(t)}|\grad f|,$$ for any $t\in (\inf f,\sup f).$\end{lem}

\section{An one end result}

\begin{theorem} Let $M$ be a complete noncompact 16-dimensional manifold with holonomy group Spin(9).
Assume that the lowest spectrum satisfies $\lambda_1(M)\geq
\frac{216}{7}$ . Then $M$ has only one end with infinite volume.
\end{theorem}
\beg{proof} Suppose that $M$ has at least two infinite volume ends,
$E_1, E_2$. Since $\lambda_1(M)>0,$ $E_1, E_2$ must be nonparabolic.
Let $f$ be the harmonic function constructed as in the previous
section. Let $e_1=\frac{\grad f}{|\grad f|}$ and $\{e_1,\cdots, e_8,
e_9,\cdots,e_{16}\}$ be a local orthonormal frame as in theorem
\ref{thm:cayley harmonic} such that $e_1f=|\grad f|,\  e_\alpha
f=0,\  2\leq \alpha\leq 16$ at a point $x$ and
$$\sum_{i=1}^8f_{ii}=0,$$ hence we have \beay \sum_{i,j=1}^{16}f_{ij}^2&\geq &f_{11}^2+\sum_{i=2}^8f_{ii}^2+2\sum _{i=2}^{16}f_{1j}^2\\
&\geq  &f_{11}^2+\frac{1}{7}\lt(\sum_{i=2}^8f_{ii}\rt)^2+2\sum
_{i=2}^{16}f_{1j}^2\\
&\geq &\frac{8}{7}\sum_{j=1}^8f_{1j}^2\\
&=&\frac{8}{7}|\grad|\grad f||^2\eeay at $x$. Combining the above
inequality with Bochner formula gives us \beay \frac{1}{2}\lap
|\grad f|^2&=&\sum _{i,j=1}^{16}f_{ij}^2+\ric(\grad f,\grad
f)\\
&\geq &\frac{8}{7}|\grad|\grad f|^2|^2-36|\grad f|^2.\eeay Let
$g=|\grad f|^{6/7},$ the above inequality becomes \bea
\label{e:1end1}\lap g\geq -\frac{216}{7}g.\eea The variational
principle of $\lambda _1(M)$ implies that for any compactly
supported smooth function $\phi \in C^\infty _c(M),$ we have
\beay \frac{216}{7}\int_M\phi ^2g^2&\leq &\int _M|\grad (\phi g)|^2\\
&=&\int_M\lt( |\grad \phi |^2 g^2+|\grad g|^2\phi
^2+\frac{1}{2}\la\grad \phi
^2,\grad g^2\ra \rt) \\
&=&\int_M|\grad \phi |^2g^2-\int_M\phi ^2g\lap g.\eeay Combining the
above with (\ref{e:1end1}), we have \bea \label{e:1end2} 0&\leq
&\int_M\phi ^2g\lt(\lap
g+\frac{216}{7}g\rt)\\
\nonum&\leq &\int_M|\grad \phi ^2|g^2.\eea We choose $\phi =\psi
\cdot \chi$ to be the product of two compactly smooth functions. For
any $\eps \in (0,1/2),$ we construct $\psi ,\chi$ as follows \beay
\chi (x)=\lt\{\bay{cll}0
&\mbox{on}& \lvl(0, \eps/2)\cup \lvl(1-\eps/2,1)\\
(\log 2)^{-1}(\log f-\log (\eps /2))&\mbox{on}&\lvl(\eps/2,\eps)\cap
(M\setminus E_1)\\
(\log 2)^{-1}(\log (1-f)-\log (\eps
/2))&\mbox{on}&\lvl(1-\eps,1-\eps/2)\cap
E_1\\
1&\mbox{otherwise}&\eay\rt..\eeay \beay \psi
=\lt\{\bay{clc}1&\mbox{on}&B(R-1)\\
R-r&\mbox{on}&B(R)\setminus B(R-1)\\
0&\mbox{on}&M\setminus B(R)\\\eay\rt..\eeay Then applying the right hand side of
(\ref{e:1end2}), we have \bea \label{e:1end3}\int_M|\grad \phi |^2g^2\leq 2\int_M|\grad
\psi|^2\chi^2|\grad f|^{\frac{12}{7}}+2\int_M|\grad \chi|^2\psi^2|\grad f|^{\frac{12}{7}}.\eea
$M$ is Einstein and the Ricci curvature satisfies $\ric_M=-36$ under our normalization. The
local gradient estimate of Cheng-Yau \cite{CY75} (see also \cite{liwangpos2}) implies
that\beay |\grad f|&\leq &C f, \eeay for some constant $C$. The above inequality implies
that\bea \label{e:1end CY}|\grad f|&\leq &C|1-f|,\eea by replacing $f$ with $1-f.$ On $E_1$,
the first term of (\ref{e:1end3}) can be estimated by \bea\label{e:1end4} \int_{E_1}|\grad
\psi|^2\chi^2|\grad f|^{\frac{12}{7}}&\leq & \lt(\int_\Omega |\grad
f|^2\rt)^{6/7}\lt(\int_\Omega 1\rt)^{1/7},\eea where $\Omega =E_1\cap (B(R)\setminus
B(R-1))\cap (\lvl(1-\eps,1-\eps/2)\cup\lvl(\eps/2,\eps)).$ Since $0<\eps<1/2,$ $\eps/2\leq
1-f$ on $\Omega $ and we have \beay \int_\Omega 1&\leq &\int_\Omega
\lt(\frac{2(1-f)}{\eps }\rt)^2\\
&\leq &C_1\eps ^{-2}\exp (-2\sqrt {\lambda _1(M)}R),\eeay where the
last inequality follows from lemma \ref{lem:est1-pos}. Combining
lemma \ref{lem:est1-pos}, the above inequality and (\ref{e:1end4}),
we conclude that \bea \label{e:1end5}\int_{E_1}|\grad \psi |^2\chi
^2|\grad f|^{12/7}\leq C_2\eps ^{-2/7}\exp (-2\sqrt {\lambda
_1(M)}R).\eea The second term of (\ref{e:1end3}) can be estimated by
\beay \int_{E_1}|\grad \chi |^2\psi^2|\grad f|^{12/7}&\leq &(\log
2)^{-2}\int_{\lvl(1-\eps, 1-\eps/2)\cap
E_1\cap B(R)}|\grad f|^{12/7}\lt|\grad \log (1-f)\rt|^2\\
&=&(\log 2)^{-2}\int_{\lvl(1-\eps, 1-\eps/2)\cap E_1\cap B(R)}|\grad
f|^{2+12/7}(1-f)^{-2}\\
&\leq &C_3\int_{\lvl(1-\eps, 1-\eps/2)\cap E_1\cap B(R)}|\grad
f|^2(1-f)^{-2/7},\eeay where the last inequality follows from
(\ref{e:1end CY}). Co-area formula and lemma \ref{lem:est-weight}
give us\beay \int_{\lvl(1-\eps, 1-\eps/2)\cap E_1\cap B(R)}|\grad
f|^2(1-f)^{-2/7}&=&\int_{1-\eps}^{1-\eps/2}(1-t)^{-2/7}\int_{l(t)\cap
E_1\cap
B(R)}|\grad f|dAdt\\
&\leq &\int_{l(b)}|\grad f|dA\int_{1-\eps}^{1-\eps/2}(1-t)^{-2/7}dt\\
&\leq &C_4\eps ^{5/7}\int_{l(b)}|\grad f|dA.\eeay Therefore, combing
the above inequalities, (\ref{e:1end3}) becomes \bea
\label{e:1end6}\int_{E_1}|\grad \phi |^2g^2\leq C_5\lt(\eps
^{-2/7}\exp (-2\sqrt {\lambda_1(M)}R)+\eps ^{5/7}\rt).\eea Applying
the same argument to $1-f$ instead of $f$ to the rest of the ends of
$M$, we have \bea \label{e:1end7}\int_{M\setminus E_1}|\grad \phi
|^2g^2\leq C_5\lt(\eps ^{-2/7}\exp (-2\sqrt {\lambda_1(M)}R)+\eps
^{5/7}\rt).\eea Combining (\ref{e:1end2}), (\ref{e:1end6}) and
(\ref{e:1end7}), letting $R\goto +\infty $ and $\eps \goto 0,$ we
conclude that $$\lap g=-\frac{216}{7}g,$$ and hence all the
inequalities in proving (\ref{e:1end1}) are indeed equalities. In
particular, $(f_{\alpha \beta})$ is diagonal and there exists a
function
$\mu $ such that \bea \label{e:1end70}(f_{\alpha \beta })=\beg{pmatrix}-7\mu&&\\
&D_1&\\
&&D_2\end{pmatrix},\eea where $D_1=\mu I$ and $D_2$ is the $8\times
8$ zero matrix. Since $f_{1\alpha }=0$ for any $\alpha \not=1,$
$|\grad f|$ is constant along the level set of $f$. In particular,
the level sets of $|\grad f|$ and $f$ coincide. Suppose $|\grad
f|(x)=0,$ by considering $f+c,$ we may assume that $f(x)=0.$ The
regularity theory of harmonic functions asserts that $f$ locally in
a neighborhood of $x$ behaves like a homogeneous harmonic polynomial
in $\re ^n$ with the origin at $x$. This is impossible since the
level sets of $|\grad f|$ and $f$ coincide. Hence $|\grad f|\not=0$
on $M$ and $M$ is diffeomorphic to $\re \times N,$ where $N$ is
given by the level set of $f$. $N$ is compact since we have assumed
that $M$ has at least two ends. Fix a level set $N$ of $f$. We
choose a local orthonormal frame $\{e_i\}_{\alpha =2}^{16}$ of $N$
and $e_1=\frac{\grad f}{|\grad f|}$. Let $\gamma (t)$ be the
integral curve of $e_1$ and $\{e_\alpha (t)\}_{\alpha =2}^{16}$ be
the parallel transport of $\{e_\alpha \}_{\alpha =2}^{16}$ along
$\gamma.$ $\la\grad _{e_1}e_1,e_1\ra=0=\la\grad _{e_1}e_1,e_\alpha
\ra$ for any $\alpha \geq 2$ implies $$\grad _{e_1}e_1=0,$$ and
hence $\gamma $ is a geodesic. The second fundamental form of the
level set of $f$ satisfies the following equations
\bea\label{e:1end8}
f_{\alpha \beta }&=&e_\alpha e_\beta f-(\grad _{e_\alpha }e_\beta )f\\
\nonum&=&\la-(\grad _{e_\alpha }e_\beta ),e_1\ra f_1\\
\nonum&=&h_{\alpha \beta }f_1 \\
\label{e:1end9}\grad _{e_\alpha }e_1&=&\sum_{\beta =2}^{16}h_{\alpha
\beta}e_\beta\eea where $h_{\alpha \beta }=\la-(\grad _{e_\alpha
}e_\beta ),e_1\ra$ is the second fundamental form of $N.$ We now
compute the curvature of $M$ \bea \label{e:1end10}\la R(e_1,e_\alpha
)e_1,e_\alpha \ra &=&\la \grad _{e_1}\grad _{e_\alpha }e_1-\grad
_{e_\alpha}\grad _{e_1
}e_1-\grad _{[e_1,e_\alpha]}e_1,e_\alpha \ra \\
\nonum&=&\la \grad _{e_1}\grad _{e_\alpha }e_1,e_\alpha \ra-\la
\grad
_{[e_1,e_\alpha]}e_1,e_\alpha \ra\\
\nonum&=&\la \grad _{e_1}\grad _{e_\alpha }e_1,e_\alpha \ra-\la
\grad _{\grad
_{e_1}e_\alpha -\grad _{e_\alpha }e_1 },e_\alpha \ra\\
\nonum&=&\la \grad _{e_1}\grad _{e_\alpha }e_1,e_\alpha
\ra+\sum_{\beta =2}^{16}\la\grad _{e_\alpha }e_1 ,e_\beta \ra
\la\grad
_{e_\beta}e_1,e_\alpha \ra \\
\nonum&=&\sum _{\beta =2}^{16}\la e_1(h_{\alpha
\beta})e_\beta,e_\alpha\ra+\sum_{\beta=2}^{16}h^2_{\alpha \beta } \\
\nonum&=&e_1(h_{\alpha \alpha })+h_{\alpha \alpha }^2,\eea where we
have used (\ref{e:1end8}), (\ref{e:1end9}) and the fact that
$(f_{\alpha \beta})$ is diagonal. In particular, since $M$ is
covered by $\cay$, (\ref{e:1end70}) and (\ref{e:1end10}) implies the
sectional curvature
$$K_M(e_1,e_k)=K_{\cay}(e_1,e_k )=0,$$ for any $k\geq 9.$ It contradicts to the fact that the sectional curvature
of $\cay$ is pinched between $-4$ and $-1.$ Therefore, $M$ has only
one infinite volume end.
\end{proof}


\section{Splitting type theorem}
\begin{theorem} \label{thm:one end} Let $M$ be a complete noncompact 16-dimensional manifold with holonomy group Spin(9).
Assume that the lowest spectrum of $M$ achieves the maximal value,
that is $\lambda_1(M)=121$. Then either

(1) $M$ has only one end; or

(2) $M$ is diffeomorphic to $\re \times N$ with metric
$$ds_M^2=dt^2+e^{-4t}\sum_{k=2}^8\omega _k^2+e^{-2t}\sum_{k=9}^{16}\omega _k^2,$$ where
$\{\omega _2,\cdots, \omega _{16}\}$ is an orthonormal basis for a
compact manifold $N$ given by a compact quotient of the horosphere
of the universal cover $\tilde M$ of $M$.
\end{theorem}
\begin{proof} Since $\lambda _1(M)>0,$ $M$ is nonparabolic and hence
$M$ has at least one nonparabolic end. Assume that $M$ has at least
two ends. Theorem \ref{thm:one end} implies that $M$ must have a
parabolic end. Let $E_1$ be a nonparabolic end and $E_2$ be a
parabolic end with respect to $B_p(R_0)$, the geodesic ball with
radius $R_0$ centered at $p$. In other words, $E_1,E_2$ are two
unbounded component of $M\setminus B_p(R_0).$ Let $\bar \gamma
:[0,+\infty )\goto M$ be a geodesic ray with $\bar \gamma (0)=p$ and
$\bar \gamma ([R_0,+\infty))\subseteq E_2,$ for some $a>0.$ Let
$\beta (x)=\lim_{t\goto \infty }(t-r(x,\bar \gamma (t)))$ be the
Busemann function with respect to $\bar \gamma .$ Theorem
\ref{thm:lap} gives us \beay \lap r(x,\bar \gamma (t))\leq 14\coth
(2r(x,\bar \gamma (t)))+8\coth r(x,\bar \gamma (t)), \eeay which
implies
$$\lap \beta \geq -22,$$ in the sense of distribution. Let $f=\exp (11\beta
),$ we compute \beay \lap f&=&11f\lap \beta +121f|\grad
\beta|^2\\
&\geq &-121f.\eeay Using the variation principle of $\lambda
_1(M)=121,$ for any $\phi \in C^\infty _c(M)$ nonnegative smooth
function with compact
support, we have \beay 121\int_M\phi ^2f^2&\leq &\int _M|\grad(\phi f)|^2\\
&=&\int_M|\grad \phi|^2f^2+\frac{1}{2}\int_M\la\grad \phi ^2,\grad
f^2\ra
+\int_M\phi ^2|\grad f|^2\\
&=&\int_M|\grad \phi|^2f^2-\int _M\phi^2f\lap f,\eeay thus \bea \label{e:61}\int_M
\phi^2f(\lap f+121f)\leq \int_M|\grad \phi |^2f^2.\eea Follow the argument in
\cite{liwangsymm}, if we choose the following cut-off function
$$\phi (x)=\left\{\bay{cl}1 &\mbox{ on  }\ \ B_p(R)\\
\frac{2R-r(x)}{R}&\mbox{ on  }\ \ B_p(2R)\setminus B_p(R)\\
0&\mbox{ on  }\ \ M\setminus B_p(R)\eay \right.,$$ then the right
hand side of (\ref{e:61}) converges to zero as $R\goto +\infty .$
Indeed, \bea \label{e:62}\int_M|\grad \phi ^2|f^2\leq
R^{-2}\int_{\lt(B_p(2R)\setminus B_p(R)\rt)\cap
E_2}f^2+R^{-2}\int_{\lt(B_p(2R)\setminus B_p(R)\rt)\setminus
E_2}f^2.\eea For an end $E$, let $V_R(E)$ be the volume of the set
$B_p(R)\cap E$ and let $k-1\leq R<k.$ The first term on the right
hand side of (\ref{e:62}) can be estimated by \bea \label{e:63}
\int_{\lt(B_p(2R)\setminus B_p(R)\rt)\cap E_2}f^2&\leq &\sum
_{i=1}^{k}\int_{(B_p(R+i)\setminus
B_p(R+i-1))\cap E_2}f^2\\
\nonum&\leq &\sum _{i=1}^{k}e^{22(R+i)}\lt(V_{E_2}(R+i)\setminus
V_{E_2}(R+i-1)\rt)\\
\nonum&\leq &C_1\sum _{i=1}^{k}e^{22(R+i)}e^{-22(R+i-1)}\\
\nonum&\leq &C_2R,\eea  where the second inequality follows from $|\beta (x)|\leq r(x,p)$ and
the third inequality follows from the volume estimate on a parabolic end $E$ of
\cite{liwangpos},
$$V_\infty(E)-V_R(E)\leq C\exp (-2\sqrt {\lambda (E)}R) \ \ \mbox{if
}\lambda(E)>0.$$ On the other hand, let $\tau $ be the geodesic ray
given in lemma \ref{lem:busemann}. For any $x\in M\setminus
(B_p(R_0)\cup E_2),$ then $\tau $ must intersect $B_p(R_0).$ Let $y$
to be the first point on $\tau $ that intersects $B_p(R_0),$
(\ref{e:bus})
implies \beay \beta (y)-\beta (x)&\geq &r(x,y)\\
&\geq &r(x,p)-r(y,p),\eeay and hence \beay
\beta (x)&\leq &-r(x,p)+r(y,p)+\beta (y)\\
\nonum &\leq &-r(x,p)+2r(y,p)\\
\nonum&\leq &-r(x,p)+2R_0,\eeay and hence the second term of the
right hand side of (\ref{e:62}) can now be estimated by
\bea\label{e:64} \int_{\lt(B_p(2R)\setminus B_p(R)\rt)\setminus
E_2}f^2&\leq &\sum
_{i=1}^{k}\int_{(B_p(R+i)\setminus B_p(R+i-1))\setminus E_2}f^2\\
\nonum&\leq &\sum
_{i=1}^{k}\int_{(B_p(R+i)\setminus B_p(R+i-1))\setminus E_2}\exp (44R_0-22r(x,p))\\
\nonum&\leq &\sum _{i=1}^{k}C_3e^{-22(R+i-1)}V(B_p(R_0+i))\\
\nonum&\leq &\sum _{i=1}^{k}C_3e^{-22(R+i-1)}e^{22(R+i)}\\
\nonum&\leq &C_4R.\eea Combining (\ref{e:62}), (\ref{e:63}) and
(\ref{e:64}), we conclude that the right hand side of (\ref{e:61})
converges to zero as $R\goto +\infty .$ Since $f$ is non-negative,
(\ref{e:61}) now implies
 \beay \lap f+121f=0, \eeay and all inequalities in the proving (\ref{e:61}) are indeed equalities and
 in particular,
  \bea \label{e:lapbeta} \lap \beta =-22,\ \ |\grad
\beta |=1,\eea and $\beta $ is smooth by the regularity of the above
equation. The above equation implies $M$ is diffeomorphic to $\re
\times N,$ where $N$ is diffeomorphic to the level set of $\beta $.
$N$ is compact since otherwise $M$ would have only one end,
contradicts to our assumption that $M$ has two ends. Let $N_0$ be
the level set of $\beta $ with $x\in N_0.$ Let $e_1=\grad \beta (x)
=\tau '(x)$ be the unit normal direction of $N_0$ at $x$, where
$\tau $ was the geodesic ray given in lemma \ref{lem:busemann}. Let
$\gamma (t)$ be the integral curve of $\grad \beta $ with $\gamma
(0)=x\in N_0$ and $e_1(t)=\gamma '(t)$. We pick a local orthonormal
frame $\{e_i\}_{A=2}^{16}$ of $N_0$ around $x$ as in the proof of
proposition \ref{p:lap} such that \beay
{R_{1i1i}}(x)&=& -4, \ 2\leq i\leq 8\\
R_{1\alpha 1\alpha }(x)&=&-1,\ 9\leq \alpha \leq 16, \eeay at $x.$
We extend the frame to a local orthonormal frame
$\{e_A(t)\}_{A=2}^{16}$ along $\gamma $ by parallel transport.
$e_1\la e_1,e_\alpha \ra =0=e_1\la e_1,e_1\ra$ implies $\grad
_{e_1}e_1=0,$ thus $\gamma (t)$ is a normal geodesic with $\gamma
'(0)=\tau '(x)$. Therefore $\gamma \equiv \tau .$  As in the proof
of proposition \ref{p:lap}, we have \beay
{R_{1i1i}}(\gamma (t))&=& -4, \ 2\leq i\leq 8\\
R_{1\alpha 1\alpha }(\gamma (t))&=&-1,\ 9\leq \alpha \leq 16, \eeay
along $\gamma .$ Bochner formula gives us \bea \label{e:hessbeta0}
0&=&\frac{1}{2}\lap |\grad
 \beta
 |^2\\
 \nonum &=&\sum _{i,j=1}^{16}\beta _{ij}^2+\ric(\grad \beta, \grad \beta
 )+\la\grad \beta ,\grad \lap \beta \ra\\
 \nonum &=&\sum _{i,j=1}^{16}\beta _{ij}^2-36.\eea The proof of proposition
 \ref{p:lap} with $e_1=\gamma '(0)=\tau '(x)$ implies \beay \beta _{11}=0,\ \ \sum _{i=2}^8\beta _{ii}&\geq& -14,\ \ \sum _{\alpha
 =9}^{16}\beta _{\alpha \alpha }\geq -8,\eeay
where the first equality comes from the fact that $\beta $ is linear
along $\tau $ (lemma \ref{lem:busemann}). Combining the above with
 (\ref{e:lapbeta}) implies \bea \label{e:hessbeta} \sum _{i=2}^8\beta _{ii}&=&
 -14,\ \  \sum _{\alpha =9}^{16}\beta _{\alpha \alpha }= -8.\eea Combining
(\ref{e:hessbeta0}) and (\ref{e:hessbeta}), we have \beay
36&=&\sum_{A,B=1}^{16}\beta _{AB}^2\\
&\geq &\sum_{i=2}^8\beta _{ii}^2+\sum_{\alpha =9}^{16}\beta^2
_{\alpha \alpha
}\\
&\geq &\frac{1}{7}\lt(\sum_{i=2}^8\beta
_{ii}\rt)^2+\frac{1}{8}\lt(\sum_{\alpha =9}^{16}\beta _{\alpha
\alpha
}\rt)^2\\
&=&36.\eeay Therefore all inequalities in the above proof are indeed
equalities $\beta _{AB}$ is diagonal and \bea
\label{e:betamatrix}\beta _{AB}=-c_A\delta _{AB},\eea where
$$c_A=\lt\{\bay{clc} 0&&A=1\\
2&&2\leq A\leq 8\\
1&&9\leq A\leq 16\eay \rt..$$ The second fundamental form of the
each level set $N_t=\{x\in M:\beta (x)=t\}$ with respect to the
normal vector $\grad \beta $ can now be
calculated \beay h_{\sigma \tau }&=&\la-\grad_{e_\sigma }e_\tau,e_1\ra\\
&=&\la-\grad_{e_\sigma }e_\tau,\grad \beta \ra\\
&=& -(\grad_{e_\sigma }e_\tau)\beta \\
&=&\beta _{\sigma \tau},\eeay where $2\leq \sigma,\tau \leq 16$ and
the last
 equality follows from the fact that $N_t$ is a level set of $\beta.$  In
 particular, we have \bea \label{e:II level set} \grad_{e_\sigma  }e_1=\sum
 _{\tau =2}^{16}\beta _{\sigma \tau}e_\tau .\eea
For any $p\in N_0,$ let $\gamma (t)$ be the integral curve of $\grad
\beta $ with $\gamma (0)=p.$ Define $\bar \psi _t(p)=\gamma (t),$
and it induces a map $\psi _t:N_0\goto N_t.$ As we have already seen
that the integral curve of $\grad \beta $ is a normal geodesic,
$\sigma (t)=\psi _t(\cdot )$ is always a normal geodesic and thus
$\psi _t$ is a geodesic flow on $M$, therefore $d\psi _t(X)$ is a
Jacobi field along each integral curve. Let $\bar e_k$ be the
restriction of $e_k$ on $N_0$, $\ 1\leq k\leq 16.$ We claim that
$d\psi_t(\bar e_i)=V_i(t),$ where
$$V_A(t)=e^{-c_At}e_A(t), \ \ 2\leq A\leq 16.$$ By the uniqueness of Jacobi
field, it is sufficient to show that $V_A(t)$ satisfies the Jacobi
equation with the same initial conditions as $d\psi _t(\bar e_A)$.
We have \beay
\grad _{\gamma '}\grad_{\gamma '} V_A&=&-c_A^2e_A\\
&=&R_{1A1A}e_A\\
&=&R(\gamma ',V)\gamma ',\eeay since $R_{AB}=R_{1A1B}$ is diagonal.
On the other hand, $ V_A(0)=\bar e_A=d\psi _0(\bar e_A)$ and
(\ref{e:II level set}) implies \beay \grad _{\gamma '}(d \psi_t(\bar
e_A))(0) &=&\grad _ {\bar e_A}
e_1(0)\\
&=&\sum _{\tau =2}^{16}\beta _{A\tau }\bar e_\tau \\
&=&-c_A\bar e_A,\eeay since we can view $e_1$ and $d \psi_t(\bar
e_A)$ as tangent vectors of a map from a rectangle. Therefore $
V'_A(0)=-c_A\bar e_A =\grad _{\gamma '}(d \psi_t(\bar e_A))(0).$ In
conclusion, each $N_t$ can be viewed as a copy of $N_0$ and $M$ is
diffeomorphic $\re \times N_0$ with metric
$$ds_M^2=dt^2+e^{-4t}\sum_{k=2}^8\omega _k^2+e^{-2t}\sum_{k=9}^{16}\omega
_k^2,$$ where $\{\omega _k\}_{k=2}^{16}$ is the coframe of $\{\bar
e_2,\cdots,\bar e_{16}\}$.

\end{proof}


   \bibliographystyle{mrl}
  \bibliography{bib}

\noindent
\address{Department of Mathematics \\ National Cheng Kung University, Taiwan \\ 1 University
road, Tainan 701, Taiwan}

\noindent \email{khlam@alumni.uci.edu}

  \end{document}